\documentclass{amsart}
\usepackage{amssymb}
\usepackage{latexsym}
\usepackage{amsthm}
\usepackage{amscd}
\theoremstyle{plain}

\newtheorem{thm}{Theorem}[section]

\newtheorem{lem}[thm]{Lemma}
\newtheorem{prop}[thm]{Proposition}

\newtheorem{defn}{Definition}[section]
\newtheorem{rem}{Remark}[section]

\numberwithin{equation}{section}

\newcommand{\md}{\mathrm{mod}}
\newcommand{\sgn}{\mathrm{sgn}}

\newcommand{\tn}{{\tilde{n}}}

\newcommand{\ZZ}{\mathbb Z}
\newcommand{\CC}{\mathbb C}

\newcommand{\RR}{\mathbb R}

\newcommand{\csp}{\null\hskip 20pt}
\newcommand{\ccsp}{\null\hskip 40pt}
\newcommand{\cccsp}{\null\hskip 80pt}
\newcommand{\ccccsp}{\null\hskip 160pt}

\allowdisplaybreaks[1]

\begin{document}

\title[HECKE OPERATORS ON WEIGHTED DEDEKIND SYMBOLS]
{Hecke operators on weighted Dedekind symbols}
\author{Shinji Fukuhara}
\thanks{The author wishes to thank Professor N.~Yui for her helpful advice.}
\thanks{\it{Address.} \rm{Shinji Fukuhara, Department of Mathematics,
  Tsuda College, Tsuda-machi 2-1-1, \\
  Kodaira-shi, Tokyo 187-8577, Japan
  (e-mail: fukuhara@tsuda.ac.jp).}
}

\begin{abstract}
Dedekind symbols generalize the classical Dedekind sums (symbols).
The symbols are determined uniquely by their reciprocity laws
up to an additive constant.
There is a natural isomorphism between the space of
Dedekind symbols with polynomial (Laurent polynomial)
reciprocity laws
and the space of cusp (modular) forms.
In this article we introduce Hecke operators on the space
of weighted Dedekind symbols.
We prove that these newly introduced operators
are compatible with Hecke operators on
the space of modular forms.
As an application, we present formulae to give Fourier
coefficients of Hecke eigenforms. In particular we give
explicit formulae for generalized Ramanujan's tau functions.
\end{abstract}

\maketitle

\section{Introduction and statement of results}
\label{sect1}

This article is a continuation of our study (\cite{F1,F2,F3})
on Dedekind symbols and modular forms.
Here we introduce and investigate Hecke operators
on Dedekind symbols, and investigate their properties.

First let us recall a few definitions in \cite{F1}
which are necessary for our subsequent discussions.
A {\em Dedekind symbol} is a generalization
of the classical Dedekind sums (\cite{RG1}),
and is defined
as a complex valued function $D$ on
  $V:=\{(p,q)\in \ZZ^+\times \ZZ\,|\,\gcd(p,q)=1\}$
satisfying
\begin{equation}\label{eqn1.1}
  D(p,q)=D(p,q+p).
\end{equation}
The symbol $D$ is determined uniquely by its {\em reciprocity law}:
\begin{equation}\label{eqn1.2}
  D(p,q)-D(q,-p)=R(p,q)
\end{equation}
up to an additive constant. The function $R$ is defined on
  $U:=\{(p,q)\in \ZZ^+\times \ZZ^+\,|\,\gcd(p,q)=1\},$
and is called a {\em reciprocity function} associated with the
Dedekind symbol $D$.
The function $R$ necessarily satisfies
the equation:
\begin{equation}\label{eqn1.3}
  R(p+q,q)+R(p,p+q)=R(p,q).
\end{equation}

When the reciprocity function $R$ is a (Laurent) polynomial in $p$ and $q$,
the symbol $D$ is called a {\em Dedekind symbol with $($Laurent$)$ polynomial
reciprocity law}. Those symbols are particularly important because
they naturally correspond to modular forms (explicit forms of such
Dedekind symbols were given in \cite{F2,F3}).

The aim of this article is to define
{\em Hecke operators on Dedekind symbols}
which are compatible with Hecke operators on modular forms.
We then apply those operators to express
Fourier coefficients of Hecke eigenforms.
For this purpose it is necessary to extend
the domain
  $V=\{(p,q)\in \ZZ^+\times \ZZ\,|\,\gcd(p,q)=1\}$
for Dedekind symbols
to $\ZZ^+\times\ZZ$.
That is, we need to define Dedekind symbol
$D(p,q)$ when $\gcd(p,q)>1$.
Thus we reach the following
definition of {\em weighted Dedekind symbols}
(hereafter we always assume an integer $w$ to be even and positive).

\begin{defn}\label{defn1.1}
A complex valued function $E$ on
$\,\ZZ^+\times\ZZ$
is called a weighted Dedekind symbol of weight $w$
if it satisfies the following two conditions:
\begin{equation}\label{eqn1.4}
  E(h,k)=E(h,k+h)
\end{equation}
for any $(h,k)\in\ZZ^+\times\ZZ$;
\begin{equation}\label{eqn1.5}
  E(ch,ck)=c^wE(h,k)
\end{equation}
for any $(h,k)\in\ZZ^+\times\ZZ$ and $c\in\ZZ^+$.

Moreover, a weighted Dedekind symbol $E$
is said to be even $($resp. odd$)$ if $E$ satisfies
\begin{equation}\label{eqn1.6}
  E(h,-k)=E(h,k) \text{\ \ \ \ $($resp.\ $E(h,-k)=-E(h,k))$}
\end{equation}
for any $(h,k)\in\ZZ^+\times\ZZ$.
\end{defn}

Roughly speaking,
a symbol $E$ is determined by its {\em reciprocity law}
\begin{equation*}
  E(h,k)-E(k,-h)=S(h,k)
\end{equation*}
up to addition of scalar multiples of
the ``trivial'' weighted Dedekind symbol.
Here $S$ is a complex valued function
defined on $\ZZ^+\times\ZZ^+$.

More precisely, let $E$ and $E'$ be Dedekind symbols of weight $w$
which have the identical reciprocity function,
namely
\begin{equation*}
  E(h,k)-E(k,-h)=E'(h,k)-E'(k,-h)
\end{equation*}
for any $(h,k)\in\ZZ^+\times\ZZ$;
then it holds that
\begin{equation*}
  E-E'=c\,G_w
\end{equation*}
where $c$ is a constant and
$G_w$ is the ``trivial'' Dedekind symbol (of weight $w$)
defined by
\begin{equation}\label{eqn1.7}
  G_w(h,k):=\left\{\gcd(h,k)\right\}^w
\end{equation}
for any $(h,k)\in\ZZ^+\times\ZZ$.

Next we would like to demonstrate the relationship between modular forms
and {\em weighted Dedekind symbols} in order to define
compatible Hecke operators
(refer to \cite{F1} for the relationship between modular forms
and {\em non-weighted Dedekind symbols}).
However, the case of non-cusp forms seems to be too involved
to treat here.
Hence we consider only the case of cusp forms in this section,
and leave
the general case to the later sections.
The statement of our results requires the following notation:
{\allowdisplaybreaks
  \begin{align*}
  \Gamma&:=SL_2(\ZZ) \text{\ (the full modular group)}, \\
  S_{w+2}&:=
    \text{the space of cusp forms on $\Gamma$ with weight $w+2$,} \\
  \mathcal{W}_w&:=\{
    W\ |\
    \text{$W$ is a Dedekind symbol of weight $w$}
    \}, \\
  \mathcal{W}_w^-&:=\left\{
    W\in\mathcal{W}_w|\ W \text{\ \ is odd\ }
    \right\}, \\
  \mathcal{W}_w^+&:=\left\{
    W\in\mathcal{W}_w|\ W \text{\ \ is even\ }
    \right\}, \\
  \mathcal{E}_w&:=\{
    E\ |\
    \text{$E$ is a Dedekind symbol of weight $w$ such
          that $E(h,k)-E(k,-h)$ is } \\
    &\cccsp\csp\text{\ \ \ a homogeneous polynomial in $h$ and $k$
      of degree $w$}
    \} \\
    &\ \ \text{(an element of $\mathcal{E}_w$ is essentially
    a period polynomial modulo $h^w-k^w$ \cite{F1,KZ1}),} \\
  \mathcal{E}_w^-&:=\left\{
    E\in\mathcal{E}_w|\ E \text{\ \ is odd\ }
    \right\}, \\
  \mathcal{E}_w^+&:=\left\{
    E\in\mathcal{E}_w|\ E \text{\ \ is even\ }
    \right\}, \\
  \mathcal{U}_w&:=\{
    g\ |\
    \text{$g$ is a homogeneous polynomial
          in $h$ and $k$ of degree $w$} \\
    &\ccsp\csp\text{satisfying $g(h+k,k)+g(h,h+k)=g(h,k)$ and $g(1,1)=0$}
    \}, \\
  \mathcal{U}_w^-&:=\left\{
    g\in \mathcal{U}_w|\
    \text{\ $g$ is an odd polynomial, i.e., $g(h,-k)=-g(h,k)$}
    \right\}, \\
  \mathcal{U}_w^+&:=\left\{
    g\in \mathcal{U}_w|\
    \text{\ $g$ is an even polynomial, i.e., $g(h,-k)=g(h,k)$}
    \right\}. \\
  \end{align*}
}
It is obvious that
$\mathcal{W}_w^+\oplus\mathcal{W}_w^-=\mathcal{W}_w$,\ \
$\mathcal{E}_w^+\oplus\mathcal{E}_w^-=\mathcal{E}_w$,\ \
$\mathcal{U}_w^+\oplus\mathcal{U}_w^-=\mathcal{U}_w$\ \
and\ \
$\mathcal{E}_w^\pm\subset\mathcal{W}_w^\pm$.

For a cusp form $f\in S_{w+2}$ and $(h,k)\in\ZZ^+\times\ZZ$,
we define $E_f$ by
\begin{equation}\label{eqn1.8}
  E_f(h,k)=\int_{k/h}^{i\infty}f(z)(hz-k)^{w}dz.
\end{equation}
Furthermore we define $E_f^-$ and $E_f^+$, respectively, by
\begin{equation*}
  E_f^-(h,k)=\frac{1}{2}\{E_f(h,k)-E_f(h,-k)\}
\end{equation*}
and
\begin{equation*}
  E_f^+(h,k)=\frac{1}{2}\{E_f(h,k)+E_f(h,-k)\}
\end{equation*}
for any $(h,k)\in \ZZ^+\times\ZZ$.
Then it is shown that $E_f$ is a Dedekind symbols of weight $w$,
and we can define maps
\begin{equation*}
  \alpha_{w+2}:S_{w+2}\to\mathcal{W}_w,\ \
  \alpha_{w+2}^\pm:S_{w+2}\to\mathcal{W}_w^\pm
\end{equation*}
by
\begin{equation*}
  \alpha_{w+2}(f)=E_f,\ \
  \alpha_{w+2}^\pm(f)=E_f^\pm.
\end{equation*}
Furthermore, we know that $E_f$ and $E_f^{\pm}$ have polynomial
reciprocity laws,
that is, $E_f\in\mathcal{E}_w$ and
$E_f^{\pm}\in\mathcal{E}_w^\pm$.
Hence we have the restricted maps
\begin{equation*}
  \alpha_{w+2}^\pm:S_{w+2}\to\mathcal{E}_w^\pm
\end{equation*}
(we use the same notation $\alpha_{w+2}^\pm$ for the restricted maps).
Using the trivial element $F_w\in \mathcal{E}_w^+$
defined (for any $(h,k)\in \ZZ^+\times\ZZ$) by
\begin{equation}\label{eqn1.9}
  F_w(h,k)=h^w,
\end{equation}
we obtain the following:

\begin{thm}\label{thm1.1}
The map
\begin{equation*}
  \alpha_{w+2}^-:S_{w+2}\to\mathcal{E}_w^-
\end{equation*}
is an isomorphism $($between vector spaces$)$ and the map
\begin{equation*}
  \alpha_{w+2}^+:S_{w+2}\to\mathcal{E}_w^+
\end{equation*}
is a monomorphism such that
the image $\alpha_{w+2}^+(S_{w+2})$
is a subspace of $\mathcal{E}_w^+$ of codimension two,
and that $\alpha_{w+2}^+(S_{w+2})$,
$F_w$ and $G_w$ span $\mathcal{E}_w^+$.
\end{thm}

Next we will see how weighted Dedekind symbols are linked
to reciprocity functions.
For a weighted Dedekind symbol $E$, let $\beta_{w}(E)$
be defined by
  \begin{equation*}
    \beta_{w}(E)(h,k)=E(h,k)-E(k,-h)
  \end{equation*}
for any $(h,k)\in \ZZ^+\times\ZZ$.
In other words, $ \beta_{w}(E)$ is the reciprocity function of $E$.
In the case of Dedekind symbol $E_f$ associated with
a cusp form $f$, this has the following expression:
\begin{equation}\label{eqn1.10}
  \beta_w(E_f)(h,k)=\int_{0}^{i\infty}f(z)(hz-k)^{w}dz.
\end{equation}
Hence obviously $\beta_w(E_f)(h,k)$ is a homogeneous polynomial
in $h$ and $k$. Furthermore we know $\beta_w(E_f)\in\mathcal{U}_w$.
Thus we have a homomorphism
  \begin{equation*}
    \beta_{w}:\mathcal{E}_w\to \mathcal{U}_w.
  \end{equation*}
Our second result is:

\begin{thm}\label{thm1.2}
The homomorphism $\beta_{w}:\mathcal{E}_w\to \mathcal{U}_w$
is an epimorphism such that
$\beta_{w}(\mathcal{E}_w^\pm)=\mathcal{U}_w^\pm$
and $\ker\beta_{w}$ is one dimensional subspace of $\mathcal{E}_w$
spanned by $G_w$.

In particular, the restricted map
  \begin{equation*}
    \beta_{w}^-:\mathcal{E}_w^-\to \mathcal{U}_w^-
  \end{equation*}
is an isomorphism, and
  \begin{equation*}
    \beta_{w}^+:\mathcal{E}_w^+\to \mathcal{U}_w^+
  \end{equation*}
is an epimorphism such that
$\ker\beta_{w}^+$ is one dimensional subspace of $\mathcal{E}_w^+$
spanned by $G_w$.
\end{thm}

Here we examine the composed maps
\begin{equation*}
  \beta_{w}^\pm\alpha_{w+2}^\pm:S_{w+2}\to
  \mathcal{E}_w^\pm\to\mathcal{U}_w^\pm.
\end{equation*}
Note that $\beta_w^\pm(E_f^\pm)(h,k)=\beta_w^\pm\alpha_{w+2}^\pm(f)(h,k)$
is a homogenized form of the period polynomial
(Kohnen-Zagier \cite[pp.\ 199--200]{KZ1}) for $f$.
This means the composed maps
\begin{equation*}
  \beta_{w}^-\alpha_{w+2}^-:S_{w+2}\to
  \mathcal{E}_w^-\to\mathcal{U}_w^-
\end{equation*}
and
\begin{equation*}
  \beta_{w}^+\alpha_{w+2}^+:S_{w+2}\to
  \mathcal{E}_w^+\to\mathcal{U}_w^+
\end{equation*}
can be identified with the Eichler-Shimura
isomorphisms (refer to \cite[p.\ 200]{KZ1}, \cite[Theorem 7.3]{F1}).
In fact,
$\beta_{w}^-\alpha_{w+2}^-$ is an isomorphism, and
$\beta_{w}^+\alpha_{w+2}^+$
is an monomorphism such that the image
$\beta_{w}^+\alpha_{w+2}^+(S_{w+2})$
and $h^w-k^w$ span $\mathcal{U}_w^+$.

These facts may be summarized in the following commutative diagram:

\setlength{\unitlength}{1mm}
\begin{picture}(125,50)(4,-1)
  \put(44,35){\framebox(39,10)
    {\shortstack{the space of cusp forms \\
      of weight $w+2$}}
  }
  \put(0,5){\framebox(65,15)
    {\shortstack{the space of odd (resp. even) Dedekind \\
                 symbols of weight $w$ with polynomial \\
                 reciprocity laws ($\md$ $F_w$ and $G_w$ if even)}}
  }
  \put(77,5){\framebox(50,15)
    {\shortstack{the space of odd (resp. even) \\
                 period polynomials of degree $w$ \\
                 ($\md$ $h^w-k^w$ if even).}}
  }
  \put(44,33){\vector(-1,-1){10}}
  \put(28,27){\makebox(10,5){$\alpha_{w+2}^\pm$}}
  \put(38,25){\makebox(10,5){$\cong$}}
  \put(84,33){\vector(1,-1){10}}
  \put(97,22){\makebox(10,15)
    {\shortstack{the Eichler-Shimura \\
                 isomorphism}}
  }
  \put(80,25){\makebox(10,5){$\cong$}}
  \put(68,12){\vector(1,0){7}}
  \put(66,14){\makebox(10,5){$\beta_{w}^\pm$}}
  \put(66,6){\makebox(10,5){$\cong$}}

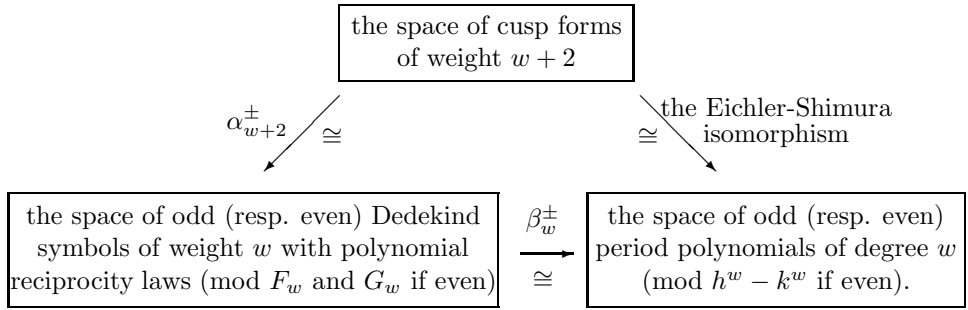
  \put(66,-3){\makebox(15,5){Diagram 1: The case of cusp forms.\ccsp }}
\end{picture}

In the above diagram, we have Hecke operators for modular forms
and period polynomials.
Indeed, Manin \cite{M1} and Zagier \cite{Z1} proved
that there are well defined Hecke operators
on period polynomials which are compatible with
the Eichler-Shimura isomorphism.
However, no such operators are yet known
for Dedekind symbols. Under these circumstances,
we introduce the following operators:

\begin{defn}\label{defn1.2}
For any positive integer $n$,
we define the operator $T_n$ on $\mathcal{W}_w$ by
\begin{equation}\label{eqn1.11}
  (T_{n}E)(h,k)
    :=\sum_{\substack{ad=n \\ 0<d}}
      \sum_{b(\md d)}E(dh,ak+bh).
\end{equation}
We may call the operator $T_n$ as the Hecke operator.
\end{defn}

We will show that $T_{n}$ maps any weighted Dedekind symbol $E$
in $\mathcal{W}_w$
onto another weighted Dedekind symbol in $\mathcal{W}_w$.
Furthermore we will show that $T_{n}$ preserve
$\mathcal{W}_w^\pm$, namely $T_{n}$ induces operators
on $\mathcal{W}_w^\pm$:
\begin{equation*}
  T_{n}:\mathcal{W}_w^\pm\to\mathcal{W}_w^\pm.
\end{equation*}
Then we have the following result which asserts that
Hecke operators on Dedekind symbols are compatible
with well known Hecke operators on cusp forms:

\begin{thm}\label{thm1.3}
The following diagram commutes:
\begin{equation*}
  \begin{CD}
  S_{w+2} @>{\alpha_{w+2}^\pm}>> \mathcal{W}_w^\pm \\
  @VV{T_n}V @VV{T_n}V \\
  S_{w+2} @>{\alpha_{w+2}^\pm}>> \mathcal{W}_w^\pm.
  \end{CD}
\end{equation*}
\end{thm}

To ease the notation, We will use the same notation $T_n$
for the Hecke operators on $S_{w+2}$ and also on $\mathcal{W}_w$.

Finally, as an application of Hecke operators on Dedekind symbols,
we present formulae giving Fourier coefficients of cusp forms
which are Hecke eigenforms
in terms of Dedekind symbols:

\begin{thm}\label{thm1.4}
Let
\begin{equation*}
  f(z)=\sum_{n=1}^{\infty}a_f(n)e^{2\pi inz}\in S_{w+2}
\end{equation*}
be a normalized Hecke eigenform having $a_f(n)$
as its $n$th Fourier coefficient.
Then there exists $(h,k)\in \ZZ^+\times\ZZ$ such that
$E_f(h,k)\ne 0$, and for such $(h,k)$, it holds that
\begin{equation}\label{eqn1.12}
  a_f(n)=\frac{T_nE_f(h,k)}{E_f(h,k)}.
\end{equation}
\end{thm}

When we have an explicit description for $E_f$,
the expression \eqref{eqn1.12} is
very useful to calculate the Fourier coefficient $a_f(n)$
for any $n\geq1$.

For example, for $f$ in
$S_{\ell+2}\ (\ell=10,14,16,18,20,24)$,
we can calculate $a_f(n)$ rather efficiently.
Here we illustrate Theorem \ref{thm1.4} in the special case,
e.g. $\ell=10$ and $f=\Delta$,
where $\Delta$ is the well known
Hecke eigenform of weight $12$ for $\Gamma$ defined by
\begin{equation*}
  \Delta(z):=e^{2\pi iz}\prod_{n=1}^{\infty}(1-e^{2\pi inz})^{24}.
\end{equation*}
Let $\tau(n)$ be the $n$th Fourier coefficient of $\Delta$,
namely
\begin{equation*}
  \Delta(z)=\sum_{n=1}^{\infty}\tau(n)e^{2\pi inz}.
\end{equation*}
Customarily $\tau(n)$ is called Ramanujan's tau function.

Using Theorem \ref{thm1.4}
together with explicit formula for $E_f$
(Theorem \ref{thm5.2}, Lemma \ref{lem6.1}),
we obtain
a surprisingly elementary formula for Ramanujan's tau function:
\begin{thm}\label{thm1.5}
Let $n$ be a positive prime integer.
Then the Ramanujan's tau function $\tau(n)$
is expressed as
\begin{align*}
  \tau(n)=&1+n^{11}
  +\frac{691}{756}(n^5-n^{11}) \\
  &\ \
    -\frac{691}{6}
    \sum_{i=0}^{n-1}
      \sum_{\substack{\left(\begin{smallmatrix}a&b\\c&d\end{smallmatrix}\right)
                      \in \Gamma/\pm 1 \\
                      ac\ne 0 \\
                      |a|,|b|,|c|,|d|\leq 2n \\
                      (i/n+b/a)(i/n+d/c)<0 \\
                      }}
      \sgn\left(\frac{i}{n}+\frac{b}{a}\right)
        (ai+bn)^{5}(ci+dn)^{5}
\end{align*}
where $\sgn(x)$ denotes $+1$ or $-1$ according
to whether $x$ is positive or negative.
\end{thm}

Throughout the article, we use the following notation
and conventions.
We assume that $w$ is an {\em even} integer
with $w\geq 2$.
For $0\leq n\leq w$, the number
$\tn$ stands for $\tn:=w-n$.
We write $\sigma_k(n)$ for the sum of the $k$th powers
of the positive divisors of $n$.
We denote by $[x]$ the greatest integer not exceeding $x\in\RR$.
We also use the notation $\sgn(x)$
for the sign of $x$, that is, $+1$, $-1$ or $0$ depending on
$x$ is positive, negative or zero, respectively.
We denote by $B_m(x)$ (resp. $B_m$)
the $m$th Bernoulli polynomial (resp. number) and
by $\bar{B}_m(x)$ the $m$th Bernoulli function:
\begin{equation*}
  \bar{B}_m(x):=-m!
    \sum^{+\infty}_{\substack{k=-\infty\\ k\ne 0}}
      \frac{e^{2\pi ikx}}{(2\pi ik)^{m}}.
\end{equation*}
It is well-known that for $0\leq x<1$, $\bar{B}_m(x)$ reduces
to $B_m(x)$.

\begin{rem}\label{rem1.1}
The Hecke operators for the classical Dedekind sums $($\cite{RG1}$)$
and the generalized Dedekind sums $($\cite{A1}$)$
were already introduced by Knopp \cite{K1}
and Parson-Rosen \cite{PR1}, respectively. The Hecke operators
on $\mathcal{W}_w$ $($Definition \ref{defn1.2}$)$ can be regarded
as natural generalizations of those operators.
\end{rem}

\section{Weighted Dedekind symbols associated with modular forms}
\label{sect2}

In the previous section, for the sake of simplicity, we restricted
our discussions to cusp forms.
However, in the subsequent sections,
we would like to deal with non-cusp modular forms as well.
In this section we investigate the relationship
among modular forms, {\em weighted Dedekind symbols}
and period polynomials.

Let $M_{w+2}$ denote the space of modular forms,
that is,
\begin{equation*}
  M_{w+2}:=
    \text{the space of modular forms on $\Gamma$ with weight $w+2$.}
\end{equation*}
We already defined the spaces
$\mathcal{U}_w$, $\mathcal{U}_w^\pm$,
$\mathcal{E}_w$, $\mathcal{E}_w^\pm$
associated with $S_{w+2}$.
Here we will define corresponding spaces
$\hat{\mathcal{U}}_w$, $\hat{\mathcal{U}}_w^\pm$,
$\hat{\mathcal{E}}_w$, $\hat{\mathcal{E}}_w^\pm$
associated with $M_{w+2}$
(they all are naturally regarded as vector spaces over $\CC$).

Now, for a polynomial $g$ in $h$ and $k$,
we define ``Laurent polynomial'' $\hat{g}$ by
\begin{equation*}
  \hat{g}(h,k):=\frac{1}{hk}\left[g(h,k)-g(1,1)\{\gcd(h,k)\}^{w+2}\right].
\end{equation*}
Using this notation, we introduce the following spaces:
{\allowdisplaybreaks
  \begin{align*}
  \hat{\mathcal{U}}_w&:=\{
    \hat{g} |\
    \text{$g$ is a homogeneous polynomial in $h$ and $k$ of degree $w+2$} \\
    &\cccsp
    \text{satisfying $hg(h+k,k)+kg(h,h+k)=(h+k)g(h,k)$} \}, \\
  \hat{\mathcal{U}}_w^-&:=\left\{
    \hat{g}\in \hat{\mathcal{U}}_w|\
    \text{\ $\hat{g}$ is odd,
          i.e., $\hat{g}(h,-k)=-\hat{g}(h,k)$}
    \right\}, \\
  \hat{\mathcal{U}}_w^+&:=\left\{
    \hat{g}\in \hat{\mathcal{U}}_w|\
    \text{\ $\hat{g}$ is even,
          i.e., $\hat{g}(h,-k)=\hat{g}(h,k)$}
    \right\}, \\
  \hat{\mathcal{E}}_w&:=\{
    E\ |\
    \text{$E$ is a Dedekind symbol such
          that $E(h,k)-E(k,-h)\in \hat{\mathcal{U}}_w$}
    \}, \\
  \hat{\mathcal{E}}_w^-&:=\left\{
    E\in\hat{\mathcal{E}}_w|\ E \text{\ \ is odd\ }
    \right\}, \\
  \hat{\mathcal{E}}_w^+&:=\left\{
    E\in\hat{\mathcal{E}}_w|\ E \text{\ \ is even\ }
    \right\}. \\
  \end{align*}
}

It is obvious that
$\hat{\mathcal{E}}_w\subset\mathcal{W}_w$,\ \
$\hat{\mathcal{E}}_w^\pm\subset\mathcal{W}_w^\pm$,\ \
$\hat{\mathcal{E}}_w^+\oplus\hat{\mathcal{E}}_w^-=\hat{\mathcal{E}}_w$\ \
and\ \
$\hat{\mathcal{U}}_w^+\oplus\hat{\mathcal{U}}_w^-=\hat{\mathcal{U}}_w$.
One can easily check that\ \
$\hat{g}\in \hat{\mathcal{U}}_w$
satisfies
\begin{equation*}
  \hat{g}(h+k,k)+\hat{g}(h,h+k)=\hat{g}(h,k)
  \text{\ \ and\ \ }
  \hat{g}(1,1)=0.
\end{equation*}
One can also show that\ \
  $\mathcal{U}_w^+=\hat{\mathcal{U}}_w^+$\ \
and that
  $\mathcal{U}_w^-$ is a codimension one subspace of $\hat{\mathcal{U}}_w^-$
as consequences of \cite[Theorem 4.1 (5) and eqn. (7.8)]{F1}.

In this setting we reformulate the relationship
between modular forms and weighted Dedekind symbols.
Throughout the article, the $n$th Fourier coefficients
of $f\in M_{w+2}$ are expressed by $a_f(n)$. Namely
\begin{equation}\label{eqn2.1}
  f(z)=\sum_{n=0}^{\infty}a_f(n)e^{2\pi inz}.
\end{equation}
Then, for a modular form $f$, $E_f$ is defined as follows:

\begin{defn}\label{defn2.1}
Let $f\in M_{w+2}$.
Firstly, for $(p,q)\in V$,
we define $E_f(p,q)$ by
\begin{equation}\label{eqn2.2}
  \begin{split}
     E_f(p,q)
           :=&\int_{z_0}^{i\infty}
             \left\{f(z)-a_f(0)\right\}(pz-q)^wdz \\
           &+\int_{q/p}^{z_0}
             \left\{f(z)-\frac{a_f(0)}{(pz-q)^{w+2}}\right\}
               (pz-q)^wdz \\
           &-a_f(0)\left\{\frac{1}{(w+1)p}(pz_0-q)^{w+1}
             +\frac{1}{p(pz_0-q)}\right\}  \\
           &\ \ (z_0\in H
               \text{ arbitrary}).
  \end{split}
\end{equation}
The right-hand side of \eqref{eqn2.2} is independent
of the choice of $z_0$ $($refer to \cite{Z2}$)$.

Secondly, for any $(h,k)\in\ZZ^+\times\ZZ$,
we define $E_f(h,k)$ by
\begin{equation}\label{eqn2.3}
     E_f(h,k):=\{\gcd(h,k)\}^wE_f(\frac{h}{\gcd(h,k)},\frac{k}{\gcd(h,k)}).
\end{equation}

More precisely,
\begin{equation}\label{eqn2.4}
  \begin{split}
     E_f(h,k)
           =&\int_{z_0}^{i\infty}
             \left\{f(z)-a_f(0)\right\}(hz-k)^wdz \\
           &+\int_{k/h}^{z_0}
             \left\{f(z)-\frac{\{\gcd(h,k)\}^{w+2}a_f(0)}
               {(hz-k)^{w+2}}\right\}
               (hz-k)^wdz \\
           &-a_f(0)\left\{\frac{1}{(w+1)h}(hz_0-k)^{w+1}
             +\frac{\{\gcd(h,k)\}^{w+2}}{h(hz_0-k)}\right\}  \\
           &\ \ (z_0\in H
               \text{ arbitrary}).
  \end{split}
\end{equation}
\end{defn}

Note that, when $f$ is a cusp form, the formula \eqref{eqn2.4}
reduces to the formula \eqref{eqn1.8}
in the previous section.

Here we will give an alternative expression for $E_f$
which is necessary for our proof of Theorem \ref{thm3.3}.
For $(h,k)\in\ZZ^+\times\ZZ$,
we introduce the function $B_f$:
\begin{equation}\label{eqn2.5}
  B_f(s;h,k):=e^{\pi is/2}h^{s-1}\int_{0}^{\infty}
             \left\{f(it+\frac{k}{h})-a_f(0)\right\}t^{s-1}dt.
\end{equation}
It is easy to see that $B_f(s;h,k)$ is well defined for $\Re (s)\gg 0$
and has a meromorphic continuation, say $B_f^*(s;h,k)$,
to the entire complex numbers.
More explicitly, we have
\begin{equation}
  \begin{split}\label{eqn2.6}
     B_f^*(s;h,k)
           =&e^{\pi is/2}h^{s-1}\int_{t_0}^{\infty}
             \left\{f(it+\frac{k}{h})-a_f(0)\right\}t^{s-1}dt \\
           &+e^{\pi is/2}h^{s-1}\int_{0}^{t_0}
             \left\{f(it+\frac{k}{h})
              -\frac{\{\gcd(h,k)\}^{w+2}a_f(0)}{(hit)^{w+2}}\right\}t^{s-1}dt \\
           &-e^{\pi is/2}h^{s-1}a_f(0)\left\{\frac{t_0^s}{s}
             +\frac{\{\gcd(h,k)\}^{w+2}t_0^{s-w-2}}{(hi)^{w+2}(s-w-2)}\right\}  \\
           &\ \ (t_0>0 \text{ arbitrary}).
  \end{split}
\end{equation}
The right hand side of \eqref{eqn2.6}
is independent of the choice of $t_0$.
It is plain that, using \eqref{eqn2.6}, $E_f(h,k)$ can be rewritten as:
\begin{equation}\label{eqn2.7}
  E_f(h,k)=B_f^*(w+1;h,k).
\end{equation}

Our first task is to show the following lemma:
\begin{lem}\label{lem2.1}
$E_f$ is a Dedekind symbol of weight $w$, that is,
$E_f\in\mathcal{W}_w$.
\end{lem}
\begin{proof}
The equation $E_f(h,k)=E_f(h,k+h)$ comes from
modularity of $f$, in particular from the formula
$f(z+1)=f(z)$.

The second equation $E_f(ch,ck)=c^wE_f(h,k)$ follows
from \eqref{eqn2.3}.
\end{proof}

Our next task is to obtain reciprocity law for $E_f$.
We introduce a ``Laurent polynomial'' $S_f$ which turn out
to be a reciprocity function for $E_f$.

\begin{defn}\label{defn2.2}
Let $f\in M_{w+2}$.
Firstly, for $(p,q)\in V$,
we define $S_f(p,q)$ by
\begin{equation}\label{eqn2.8}
  \begin{split}
     S_f(p,q)
           :=&\int_{z_0}^{i\infty}
             \left\{f(z)-a_f(0)\right\}(pz-q)^wdz
            +\int_{0}^{z_0}
             \left\{f(z)-\frac{a_f(0)}{z^{w+2}}\right\}
               (pz-q)^wdz \\
           &-a_f(0)\left\{
             \frac{1}{(w+1)p}(pz_0-q)^{w+1}
             -\frac{1}{(w+1)q}(p-\frac{q}{z_0})^{w+1}
             \right\}
             +a_f(0)\frac{1}{pq}  \\
           &\ \ (z_0\in H
               \text{ arbitrary}).
  \end{split}
\end{equation}
Again the right-hand side of \eqref{eqn2.8} is independent
of the choice of $z_0$.

Secondly, for any $(h,k)\in\ZZ^+\times\ZZ$,
we define $S_f(h,k)$ by
\begin{equation}\label{eqn2.9}
     S_f(h,k):=\{\gcd(h,k)\}^wS_f(\frac{h}{\gcd(h,k)},\frac{k}{\gcd(h,k)}).
\end{equation}

In other words,
\begin{equation}\label{eqn2.10}
  \begin{split}
     S_f(h,k)
           =&\int_{z_0}^{i\infty}
             \left\{f(z)-a_f(0)\right\}(hz-k)^wdz
            +\int_{0}^{z_0}
             \left\{f(z)-\frac{a_f(0)}{z^{w+2}}\right\}
               (hz-k)^wdz \\
           &-a_f(0)\left\{
             \frac{1}{(w+1)h}(hz_0-k)^{w+1}
             -\frac{1}{(w+1)k}(h-\frac{k}{z_0})^{w+1}
             \right\} \\
             &+a_f(0)\frac{\{\gcd(h,k)\}^{w+2}}{hk}
             \ \ \ \ (z_0\in H \text{ arbitrary}).
  \end{split}
\end{equation}
\end{defn}

Again note that, when $f$ is a cusp form,
the formula \eqref{eqn2.10}
reduces to the formula \eqref{eqn1.10}.

Now we obtain the following reciprocity law for $E_f$:
\begin{prop}\label{prop2.2}
Let $f$ be a modular form of weight $w+2$,
that is, $f\in M_{w+2}$.
Then
\begin{enumerate}
  \item
it holds that
\begin{equation}\label{eqn2.11}
  E_f(h,k)-E_f(k,-h)=S_f(h,k);
\end{equation}
  \item
$E_f$ belongs to $\hat{\mathcal{E}}_w$.
\end{enumerate}
\end{prop}
\begin{proof}
First we express $E_f(k,-h)$ as follows:
{\allowdisplaybreaks
\begin{align*}
     E_f(k,-h)
           =&\int_{\frac{-1}{z_0}}^{i\infty}
             \left\{f(z)-a_f(0)\right\}(kz+h)^wdz \\
           &+\int_{-h/k}^{\frac{-1}{z_0}}
             \left\{f(z)-\frac{\{\gcd(h,k)\}^{w+2}a_f(0)}
               {(kz+h)^{w+2}}\right\}
               (kz+h)^wdz \\
           &-a_f(0)\left\{\frac{1}{(w+1)k}
             \left(\frac{-k}{z_0}+h\right)^{w+1}
             +\frac{\{\gcd(h,k)\}^{w+2}}{k(\frac{-k}{z_0}+h)}\right\}  \\
           =&\int_{z_0}^{0}
             \left\{f(\frac{-1}{z})-a_f(0)\right\}
               (\frac{-k}{z}+h)^w\frac{dz}{z^2} \\
           &+\int_{k/h}^{z_0}
             \left\{f(\frac{-1}{z})-\frac{\{\gcd(h,k)\}^{w+2}a_f(0)}
               {(\frac{-k}{z}+h)^{w+2}}\right\}
               (\frac{-k}{z}+h)^w\frac{dz}{z^2} \\
           &-a_f(0)\left\{\frac{1}{(w+1)k}
             \left(\frac{-k}{z_0}+h\right)^{w+1}
             +\frac{\{\gcd(h,k)\}^{w+2}}{k(\frac{-k}{z_0}+h)}\right\}  \\
           &\ \ \text{(substituting $-1/z$ for $z$)} \\
           =&\int_{z_0}^{0}
             \left\{f(z)-\frac{a_f(0)}{z^{w+2}}\right\}
               (hz-k)^wdz \\
           &+\int_{k/h}^{z_0}
             \left\{f(z)-\frac{\{\gcd(h,k)\}^{w+2}a_f(0)}
               {(hz-k)^{w+2}}\right\}
               (hz-k)^wdz \\
           &-a_f(0)\left\{\frac{1}{(w+1)k}
             \left(\frac{-k}{z_0}+h\right)^{w+1}
             +\frac{\{\gcd(h,k)\}^{w+2}}{k(\frac{-k}{z_0}+h)}\right\}  \\
           &\ \ \text{(applying the formula $f(-1/z)=f(z)z^{w+2}$).} \\
\end{align*}
}
From this, we have
{\allowdisplaybreaks
\begin{align*}
     E_f(h,k)-E_f(k,-h)
           =&\int_{z_0}^{i\infty}
             \left\{f(z)-a_f(0)\right\}(hz-k)^wdz \\
           &+\int_{k/h}^{z_0}
             \left\{f(z)-\frac{\{\gcd(h,k)\}^{w+2}a_f(0)}
               {(hz-k)^{w+2}}\right\}
               (hz-k)^wdz \\
           &-a_f(0)\left\{\frac{1}{(w+1)h}(hz_0-k)^{w+1}
             +\frac{\{\gcd(h,k)\}^{w+2}}{h(hz_0-k)}\right\}  \\
           &-\int_{z_0}^{0}
             \left\{f(z)-\frac{a_f(0)}{z^{w+2}}\right\}
               (hz-k)^wdz \\
           &-\int_{k/h}^{z_0}
             \left\{f(z)-\frac{\{\gcd(h,k)\}^{w+2}a_f(0)}
               {(hz-k)^{w+2}}\right\}
               (hz-k)^wdz \\
           &+a_f(0)\left\{\frac{1}{(w+1)k}
             \left(\frac{-k}{z_0}+h\right)^{w+1}
             +\frac{\{\gcd(h,k)\}^{w+2}}{k(\frac{-k}{z_0}+h)}\right\}  \\
           =&\int_{z_0}^{i\infty}
             \left\{f(z)-a_f(0)\right\}(hz-k)^wdz \\
            &+\int_{0}^{z_0}
             \left\{f(v)-\frac{a_f(0)}{v^{w+2}}\right\}
               (hv-k)^wdv \\
           &-a_f(0)\left\{\frac{1}{(w+1)h}(hz_0-k)^{w+1}
             -\frac{1}{(w+1)k}\left(\frac{-k}{z_0}+h\right)^{w+1}\right\} \\
           &+a_f(0)\frac{\{\gcd(h,k)\}^{w+2}}{hk} \\
           =&\ S_f(h,k).
\end{align*}
}
This proves the assertion (1).

Furthermore an easy exercise yields $S_f\in \hat{\mathcal{U}}_w$,
and this implies $E_f\in \hat{\mathcal{E}}_w$.
This proves the assertion (2).
\end{proof}

For $E_f$, $E_f^-$ and $E_f^+$ are defined as before,
namely
\begin{equation*}
  E_f^\pm(h,k)=\frac{1}{2}\{E_f(h,k)\pm E_f(h,-k)\}.
\end{equation*}
Similarly,
for $S_f$, $S_f^-$ and $S_f^+$ are defined by
\begin{equation*}
  S_f^\pm(h,k)=\frac{1}{2}\{S_f(h,k)\pm S_f(h,-k)\},
\end{equation*}
and they satisfy,
by virtue of Proposition \ref{prop2.2},
the following identities
\begin{equation}\label{eqn2.12}
  E_f^\pm(h,k)-E_f^\pm(k,-h)=S_f^\pm(h,k).
\end{equation}

Then it follows that $E_f^\pm\in\hat{\mathcal{E}}_w^\pm$
for $f\in M_{w+2}$,
and we arrive at the following definition:

\begin{defn}\label{defn2.3}
\begin{enumerate}
\item
The map
$\hat{\alpha}_{w+2}:M_{w+2}\to\mathcal{W}_w$
is defined by
\begin{equation*}
  \hat{\alpha}_{w+2}(f)=E_f;
\end{equation*}
\item
The maps
$\hat{\alpha}_{w+2}^\pm:M_{w+2}\to\mathcal{W}_w^\pm$
are defined by
\begin{equation*}
  \hat{\alpha}_{w+2}^\pm(f)=E_f^\pm.
\end{equation*}
\end{enumerate}
Since
$\hat{\alpha}_{w+2}^\pm(f)=E_f^\pm\in\hat{\mathcal{E}}_w^\pm$,
we have the restricted maps of $\hat{\alpha}_{w+2}^\pm$
from $M_{w+2}$ to $\hat{\mathcal{E}}_w^\pm$. These maps
will also be denoted by $\hat{\alpha}_{w+2}^\pm$.
\end{defn}

The following theorem asserts that
the maps $\hat{\alpha}_{w+2}^\pm$
are almost bijective.
\begin{thm}\label{thm2.3}
The map
\begin{equation*}
  \hat{\alpha}_{w+2}^-:M_{w+2}\to\hat{\mathcal{E}}_w^-
\end{equation*}
is an isomorphism $($between vector spaces$)$,
and the map
\begin{equation*}
  \hat{\alpha}_{w+2}^+:M_{w+2}\to\hat{\mathcal{E}}_w^+
\end{equation*}
is a monomorphism such that
the image $\hat{\alpha}_{w+2}^+(M_{w+2})$
is a codimension two subspace of $\hat{\mathcal{E}}_w^+$,
and that $\hat{\alpha}_{w+2}^+(M_{w+2})$ and $G_w$ span  $\hat{\mathcal{E}}_w^+$.
\end{thm}
The proof of this theorem will be given after
the proof of Theorem \ref{thm2.4}.

Next we define maps from
$\hat{\mathcal{E}}_w^\pm$
to $\hat{\mathcal{U}}_w^\pm$.
For this, we will make the following definition:
\begin{defn}\label{defn2.4}
For Dedekind symbol $E\in\hat{\mathcal{E}}_w^\pm$,
the maps $\hat{\beta}_w^\pm$ are
defined by
\begin{equation*}
  \hat{\beta}_w^\pm(E)(h,k):=E(h,k)-E(k,-h).
\end{equation*}
Then obviously we have
$\hat{\beta}_w^\pm(E)\in \hat{\mathcal{U}}_w^\pm$,
and this gives the homomorphisms
\begin{equation*}
  \hat{\beta}_w^\pm:\hat{\mathcal{E}}_w^\pm\to\hat{\mathcal{U}}_w^\pm.
\end{equation*}
\end{defn}

The following theorem asserts that the maps $\hat{\beta}_w^\pm$
are bijective modulo $G_w$:

\begin{thm}\label{thm2.4}
The map
\begin{equation*}
  \hat{\beta}_w^-:\hat{\mathcal{E}}_w^-\to\hat{\mathcal{U}}_w^-
\end{equation*}
is an isomorphism $($between vector spaces$)$ and the map
\begin{equation*}
  \hat{\beta}_w^+:\hat{\mathcal{E}}_w^+\to\hat{\mathcal{U}}_w^+
\end{equation*}
is an epimorphism such that $\ker(\hat{\beta}_{w+2}^+)$
is one-dimensional subspace of $\hat{\mathcal{E}}_w^+$
spanned by $G_w$.
\end{thm}

To establish this theorem, we first prove the following lemma:
\begin{lem}\label{lem2.5}
Let $E$ and $E'$ be Dedekind symbols of weight $w$
satisfying
\begin{equation*}
  E(h,k)-E(k,-h)=E'(h,k)-E'(k,-h)
\end{equation*}
for any $(h,k)\in\ZZ^+\times\ZZ$.
Then it holds that
\begin{equation*}
  E-E'=c\,G_w
\end{equation*}
for a constant $c\in\CC$.
\end{lem}

\begin{proof}
Let $D:V\to\CC$ and $D':V\to\CC$ be the restricted maps of
$E:\ZZ^+\times\ZZ\to\CC$ and $E':\ZZ^+\times\ZZ\to\CC$,
respectively.
Then $D$ and $D'$ are (non-weighted) Dedekind symbols
which have an identical reciprocity function.
By \cite[Theorem 5.1]{F1}, we know
\begin{equation*}
  D-D'=c
\end{equation*}
for a constant $c\in\CC$
(the appearance of a constant $c$ stems from the fact that
we do not assume $D(1,0)=0$ nor $D'(1,0)=0$).
Now we have
\begin{align*}
  E(h,k)-E'(h,k)
    =&\{\gcd(h,k)\}^wD(\frac{h}{\gcd(h,k)},\frac{k}{\gcd(h,k)}) \\
      &-\{\gcd(h,k)\}^wD'(\frac{h}{\gcd(h,k)},\frac{k}{\gcd(h,k)}) \\
    =&\{\gcd(h,k)\}^wc \\
    =&c\,G_w(h,k).
\end{align*}
This completes the proof.
\end{proof}

Now we are ready to give a proof of Theorem \ref{thm2.4}.
\begin{proof}[Proof of Theorem \ref{thm2.4}]
First we show that $\hat{\beta}_w^\pm$ are epimorphisms.
Let  $S\in\hat{\mathcal{U}}_w^\pm$.
Then, by definition,
$S(h,k)$ is expressed as
\begin{equation*}
  S(h,k)=\frac{1}{hk}\left[g(h,k)-g(1,1)\{\gcd(h,k)\}^{w+2}\right]
\end{equation*}
with a homogeneous polynomial $g$ of degree $w+2$
satisfying
\begin{equation*}
  hg(h+k,k)+kg(h,h+k)=(h+k)g(h,k).
\end{equation*}
Then it holds that
\begin{equation*}
  S(h+k,k)+S(h,h+k)=S(h,k).
\end{equation*}
By Theorem 5.1 in \cite{F1}, there is
a (non-weighted) Dedekind symbol
$D$ which satisfies $D(p,q)-D(q,-p)=S(p,q)$ for any $(p,q)\in V$.
Then we define $E$ by
\begin{equation*}
  E(h,k)
    =\{\gcd(h,k)\}^wD(\frac{h}{\gcd(h,k)},\frac{k}{\gcd(h,k)})
\end{equation*}
for any $(h,k)\in\ZZ^+\times\ZZ$.
One can easily check that $E\in \hat{\mathcal{E}}_w^\pm$ and
$\hat{\beta}_w^\pm(E)=S$.
This implies $\hat{\beta}_w^\pm$ are epimorphisms.

Lemma \ref{lem2.5} shows that $\hat{\beta}_w^-$
is an monomorphism and that $\ker(\hat{\beta}_{w+2}^+)$
is one-dimensional subspace of $\hat{\mathcal{E}}_w^+$
spanned by $G_w$.
This completes the proof.
\end{proof}

Next we give a proof of Theorem \ref{thm2.3}.
\begin{proof}[Proof of Theorem \ref{thm2.3}]
We see that the composed maps
\begin{equation*}
  \hat{\beta}_w^-\ \hat{\alpha}_{w+2}^-:
    M_{w+2}\to\hat{\mathcal{E}}_w^-\to\hat{\mathcal{U}}_w^-
\end{equation*}
and
\begin{equation*}
  \hat{\beta}_w^+\ \hat{\alpha}_{w+2}^+:
    M_{w+2}\to\hat{\mathcal{E}}_w^+\to\hat{\mathcal{U}}_w^+
\end{equation*}
are nothing but the Eichler-Shimura
isomorphisms (\cite[p.\ 200]{KZ1}, \cite[Theorem 7.3]{F1}).
Hence $\hat{\beta}_w^\pm\ \hat{\alpha}_{w+2}^\pm$ are
isomorphisms.
This implies
that $\hat{\alpha}_{w+2}^+$ is a monomorphism.
Furthermore, by Theorem \ref{thm2.4}, we know that
that $\hat{\alpha}_{w+2}^+(M_{w+2})$ and $G_w$ span
$\hat{\mathcal{E}}_w^+$.

Finally, since $\hat{\beta}_w^-$ is an isomorphism
by Theorem \ref{thm2.4}, we know that $\hat{\alpha}_{w+2}^-$ is
also an isomorphism. This completes the proof.
\end{proof}

Finally we summarize these facts
in the following commutative diagram (compare with Diagram 1):

\setlength{\unitlength}{1mm}
\begin{picture}(125,50)(4,-2)
  \put(44,35){\framebox(39,10)
    {\shortstack{the space of modular \\
      forms of weight $w+2$}}
  }
  \put(0,5){\framebox(62,15)
    {\shortstack{the space of odd (resp. even) Dedekind \\
                 symbols of weight $w$ with ``Laurent \\
                 polynomial'' reciprocity laws ($\md\ G_w$)}}
  }
  \put(77,5){\framebox(50,15)
    {\shortstack{the space of odd (resp. even) \\
                 period ``Laurent polynomials'' \\
                 of degree $w$.}}
  }
  \put(44,33){\vector(-1,-1){10}}
  \put(28,27){\makebox(10,5){$\hat{\alpha}_{w+2}^\pm$}}
  \put(38,25){\makebox(10,5){$\cong$}}
  \put(84,33){\vector(1,-1){10}}
  \put(97,22){\makebox(10,15)
    {\shortstack{the Eichler-Shimura \\
                 isomorphism}}
  }
  \put(80,25){\makebox(10,5){$\cong$}}
  \put(64,12){\vector(1,0){10}}
  \put(65,14){\makebox(10,5){$\hat{\beta}_{w}^\pm$}}
  \put(65,6){\makebox(10,5){$\cong$}}

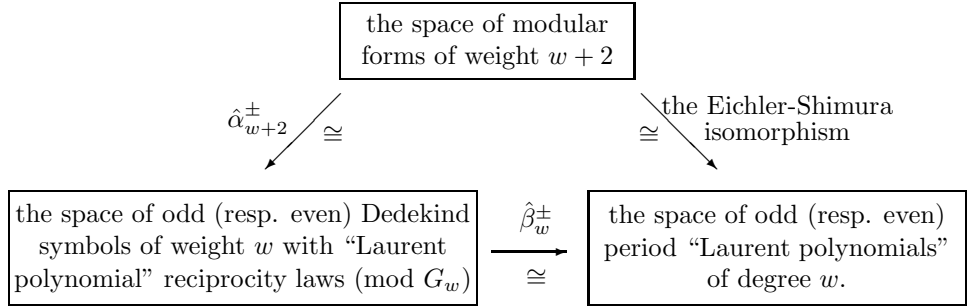
  \put(66,-3){\makebox(15,5){Diagram 2: The case of modular forms.\ccsp }}
\end{picture}

\section{Hecke operators on weighted Dedekind symbols}
\label{sect3}

One of the most important features of modular forms is
that they have Hecke operators.
The Hecke operators $T_n$ on modular forms are
defined as follows (see for example \cite{A2,S1}):

\begin{defn}\label{defn3.1}
For any $n=1,2,\ldots$,
the Hecke operator $T_n$ is defined on $M_{w+2}$ by the equation
\begin{equation*}
  (T_{n}f)(z)
    :=n^{w+1}\sum_{d|n}d^{-w-2}
      \sum_{b=0}^{d-1}f(\frac{nz+bd}{d^2}).
\end{equation*}
This can be rewritten as:
\begin{equation}\label{eqn3.1}
  (T_{n}f)(z)
    =\frac{1}{n}
    \sum_{\substack{ad=n \\ 0<d}}a^{w+2}
      \sum_{b(\md d)}f(\frac{az+b}{d}).
\end{equation}
\end{defn}

The operator $T_n$ maps the vector space $M_{w+2}$ onto itself.
In view of Diagram2,
we would like to define
Hecke operators on
the space $\mathcal{W}_w$ of Dedekind symbols with weight $w$,
compatible with Hecke operators on modular forms.

\begin{defn}\label{defn3.2}
For a positive integer $n$,
the operator $T_n$ on $\mathcal{W}_w$ is defined by the equation
\begin{equation}\label{eqn3.2}
  (T_{n}E)(h,k)
    :=\sum_{\substack{ad=n \\ 0<d}}
      \sum_{b(\md d)}E(dh,ak+bh)
\end{equation}
for any $(h,k)\in \ZZ^+\times\ZZ$.
We will call $T_n$ the Hecke the Hecke operator on weighted
Dedekind symbols.
\end{defn}

The fact that $T_{n}E$ is again a Dedekind symbol of weight $w$
will be shown by the following Lemma \ref{lem3.1}.

\begin{lem}\label{lem3.1}
Let $E$ be Dedekind symbol of weight $w$.
\begin{enumerate}
\item
Let $a$ and $d$ be fixed positive integers.
We define $\tilde{E}_{a,d}$ by
\begin{equation*}
  \tilde{E}_{a,d}(h,k)=\sum_{b(\md d)}E(dh,ak+bh)
\end{equation*}
for any $(h,k)\in \ZZ^+\times\ZZ$.
Then $\tilde{E}_{a,d}$ is also a Dedekind symbol of weight $w$.
\item
Let $n$ be a positive integer.
Then $T_{n}E$ is also a Dedekind symbol of weight $w$.
\end{enumerate}
\end{lem}

\begin{proof}
For any $(h,k)\in \ZZ^+\times\ZZ$, we see
{\allowdisplaybreaks
\begin{align*}
  \tilde{E}_{a,d}(h,k+h)
    &=\sum_{b(\md d)}E(dh,a(k+h)+bh)
    =\sum_{b(\md d)}E(dh,ak+(a+b)h) \\
    &=\sum_{b(\md d)}E(dh,ak+bh)
    =\tilde{E}_{a,d}(h,k).
\end{align*}
}
Furthermore, we have
{\allowdisplaybreaks
\begin{align*}
  \tilde{E}_{a,d}(ch,ck)
    &=\sum_{b(\md d)}E(dch,ack+bch) \\
    &=c^w\sum_{b(\md d)}E(dh,ak+bh)
    =c^w\tilde{E}_{a,d}(h,k).
\end{align*}
}
These yield the assertion (1).

The assertion (2) follows directly from (1) and
the identity
\begin{equation*}
  T_{n}E
    =\sum_{\substack{ad=n \\ 0<d}}
      \tilde{E}_{a,d}.
\end{equation*}
\end{proof}

By Lemma \ref{lem3.1},
we have just proved that $T_{n}$ is a well defined operator
on $\mathcal{W}_w$ :
\begin{equation*}
  T_{n}:\mathcal{W}_w\to\mathcal{W}_w.
\end{equation*}
Here we show that $T_{n}$ preserves $\mathcal{W}_w^\pm$:

\begin{lem}\label{lem3.2}
The Hecke operator $T_n$ preserves $\mathcal{W}_w^\pm$,
namely it holds that
\begin{equation*}
  T_n(\mathcal{W}_w^\pm)\subset\mathcal{W}_w^\pm.
\end{equation*}
\end{lem}
\begin{proof}
Let $E^\pm\in\mathcal{W}_w^\pm$.
Then we have
{\allowdisplaybreaks
\begin{align*}
  (T_nE^\pm)(h,-k)
     &=\sum_{\substack{ad=n \\ 0<d}}
            \sum_{b(\md d)}E^\pm(dh,-ak+bh) \\
     &=\sum_{\substack{ad=n \\ 0<d}}
            \sum_{b(\md d)}E^\pm(dh,-ak-bh) \\
     &=\pm\sum_{\substack{ad=n \\ 0<d}}
            \sum_{b(\md d)}E^\pm(dh,ak+bh) \\
     &=\pm (T_nE^\pm)(h,k).
\end{align*}
}
This implies $T_nE^\pm\in\mathcal{W}_w^\pm$
completing the proof.
\end{proof}

Now we formulate our theorem which asserts that Hecke
operators on Dedekind symbols are compatible with
Hecke operators on modular forms.
\begin{thm}\label{thm3.3}
Let $f$ be a modular form of weight $w+2$, that is, $f\in M_{w+2}$.
Then
\begin{enumerate}
\item
it holds that
\begin{equation*}
  \hat{\alpha}_{w+2}(T_nf)=T_n\hat{\alpha}_{w+2}(f).
\end{equation*}
In other words, the following diagram commutes
\begin{equation*}
  \begin{CD}
  M_{w+2} @>{\hat{\alpha}_{w+2}}>> \mathcal{W}_w \\
  @VV{T_n}V @VV{T_n}V \\
  M_{w+2} @>{\hat{\alpha}_{w+2}}>> \mathcal{W}_w;
  \end{CD}
\end{equation*}
\item
the Hecke operator $T_{n}$ preserves $\mathcal{W}_w^\pm$,
and the following diagram commutes
\begin{equation*}
  \begin{CD}
  M_{w+2} @>{\hat{\alpha}_{w+2}^\pm}>> \mathcal{W}_w^\pm \\
  @VV{T_n}V @VV{T_n}V \\
  M_{w+2} @>{\hat{\alpha}_{w+2}^\pm}>> \mathcal{W}_w^\pm. \\
  \end{CD}
\end{equation*}
\end{enumerate}
\end{thm}

\begin{proof}
Let $f\in M_{w+2}$. We recall the definition
of the Hecke operator $T_n$ on $f$:
\begin{equation*}
  (T_{n}f)(z)
    =\frac{1}{n}
    \sum_{\substack{ad=n \\ 0<d}}a^{w+2}
      \sum_{b(\md d)}f(\frac{az+b}{d}),
\end{equation*}
and the alternative expression \eqref{eqn2.7} for $E_{T_nf}$:
\begin{equation*}
  E_{T_nf}(h,k)=B_{T_nf}^*(w+1;h,k).
\end{equation*}

Now we have
\begin{align*}
  B_{T_nf}(s;h,k)
  &=e^{\pi is/2}h^{s-1}\int_{0}^{\infty}
    \left\{T_nf(it+\frac{k}{h})-a_{T_nf}(0)\right\}t^{s-1}dt \\
  &=e^{\pi is/2}h^{s-1}\int_{0}^{\infty}
      \sum_{\substack{ad=n \\ d>0}}
      \sum_{b(\md d)}\frac{a^{w+2}}{n}
    \left\{f(\frac{a}{d}it+\frac{ak+bh}{dh})-a_{f}(0)\right\}t^{s-1}dt \\
  &=\sum_{\substack{ad=n \\ d>0}}
      \sum_{b(\md d)}a^{w+1-s}
    e^{\pi is/2}(dh)^{s-1}\int_{0}^{\infty}
    \left\{f(ix+\frac{ak+bh}{dh})-a_{f}(0)\right\}x^{s-1}dx \\
  &\ \ \ \ \text{(substituting $at/d$ for $x$)} \\
  &=\sum_{\substack{ad=n \\ d>0}}
      \sum_{b(\md d)}a^{w+1-s}
    B_{f}(s;dh,ak+bh).
\end{align*}
Next considering the meromorphic continuations
$B_{T_nf}^*(s;h,k)$ and $B_{f}^*(s;dh,ak+bh)$ of
$B_{T_nf}(s;h,k)$ and $B_{f}(s;dh,ak+bh)$ respectively
to the entire complex numbers,
and then taking limits as $s\to w+1$,
we have
{\allowdisplaybreaks
\begin{align*}
  E_{T_nf}(h,k)
    &=B_{T_nf}^*(w+1;h,k) \\
    &=\sum_{\substack{ad=n \\ d>0}}
      \sum_{b(\md d)}
      B_{f}^*(w+1;dh,ak+bh) \\
    &=\sum_{\substack{ad=n \\ d>0}}
      \sum_{b(\md d)}
    E_{f}(dh,ak+bh) \\
    &={T_nE_f}(h,k).
\end{align*}
}
This implies that $\hat{\alpha}_{w+2}(T_nf)=T_n\hat{\alpha}_{w+2}(f)$
which proves the assertion (1).

The assertion (2) follows directly from Lemma \ref{lem3.2} and (1).
\end{proof}

\section{An application of Hecke operators on Dedekind symbols}
\label{sect4}

As an application of Hecke operators on Dedekind symbols, we present
formula which gives Fourier coefficients of Hecke eigenforms in terms of
Dedekind symbols.

\begin{thm}\label{thm4.1}
Let
\begin{equation*}
  f(z)=\sum_{n=0}^{\infty}a_f(n)e^{2\pi inz}\in M_{w+2}
\end{equation*}
be a normalized Hecke eigenform, and let $n$ be a positive integer.

Then
\begin{enumerate}
\item
  it holds that
  \begin{equation}\label{eqn4.1}
    T_nE_f^\pm(h,k)=a_f(n)E_f^\pm(h,k),
  \end{equation}
  in other words,
  \begin{equation*}
    \sum_{\substack{ad=n \\ 0<d}}
      \sum_{b(\md d)}E_f^\pm(dh,ak+bh)
      =a_f(n)E_f^\pm(h,k)
  \end{equation*}
  for any $(h,k)\in \ZZ^+\times\ZZ$;
\item
  there exists $(h,k)\in \ZZ^+\times\ZZ$ such that
  $E_f^\pm(h,k)\ne 0$.  For such $(h,k)$, it holds that
  \begin{equation*}
    a_f(n)=\frac{T_nE_f^\pm(h,k)}{E_f^\pm(h,k)}.
  \end{equation*}
\end{enumerate}
\end{thm}
\begin{proof}
Since $f$ is a normalized Hecke eigenform, it follows directly
from the definition of eigenform that
\begin{equation*}
  T_nf=a_f(n)f.
\end{equation*}
Thus we have
\begin{equation}\label{eqn4.2}
  E_{T_nf}^\pm=a_f(n)E_f^\pm.
\end{equation}
By Theorem \ref{thm3.3}, we know that
\begin{equation}\label{eqn4.3}
  E_{T_nf}^\pm=\hat{\alpha}_{w+2}^\pm(T_nf)
  =T_n\hat{\alpha}_{w+2}^\pm(f)
  =T_nE_f^\pm.
\end{equation}
From \eqref{eqn4.2} and \eqref{eqn4.3}, we have
\begin{equation}\label{eqn4.4}
  T_nE_f^\pm(h,k)=a_f(n)E_f^\pm(h,k)
\end{equation}
for any $(h,k)\in \ZZ^+\times\ZZ$.
This proves the assertion (1).

Next recall that $\hat{\alpha}_{w+2}^\pm$ are defined by
\begin{equation*}
  \hat{\alpha}_{w+2}^\pm(f)=E_f^\pm
\end{equation*}
for $f\in M_{w+2}$, and that
\begin{equation*}
  \hat{\alpha}_{w+2}^\pm :
   M_{w+2} \to \hat{\mathcal{E}}_w^\pm
\end{equation*}
are monomorphisms.
If $f$ is an eigenform, then, in particular, $f$ is non-trivial.
Hence $E_f^\pm$ are also non-trivial since $\hat{\alpha}_{w+2}^\pm$
are monomorphisms. This implies $E_f^\pm(h,k)\ne 0$
for some $(h,k)\in \ZZ^+\times\ZZ$.

When $E_f^\pm(h,k)\ne 0$, from \eqref{eqn4.1}, we have
\begin{equation*}
  a_f(n)=\frac{T_nE_f^\pm(h,k)}{E_f^\pm(h,k)}.
\end{equation*}
This proves the assertion (2) completing the proof.
\end{proof}

\section{Weighted Dedekind symbols with polynomial reciprocity laws}
\label{sect5}

In this section we give explicit description for weighted
Dedekind symbols with polynomial reciprocity laws.
Most of the arguments here are parallel to
that of the non-weighted case
so that the reader should refer to
\cite{F3,KZ1} for more details.

\begin{defn}\label{defn5.1}
Let $n$ be an integer such that $0<n<w$.
We define a sum $I_{w,n}:\ZZ^+\times\ZZ\to \CC$ by
\begin{equation*}
  I_{w,n}(h,k)
    :=
      \sum_{\substack{\left(\begin{smallmatrix}a&b\\c&d\end{smallmatrix}\right)
                      \in \Gamma/\pm 1 \\
                      ac\ne 0 \\
                      (k/h+b/a)(k/h+d/c)<0 }}
      \sgn\left(\frac{k}{h}+\frac{b}{a}\right)
        (ak+bh)^{\tn}(ck+dh)^{n}.
\end{equation*}
This sum reduces to the following finite sum
$($\cite{F3}$)$
\begin{equation*}
  I_{w,n}(h,k)
    =
      \sum_{\substack{\left(\begin{smallmatrix}a&b\\c&d\end{smallmatrix}\right)
                      \in \Gamma/\pm 1 \\
                      ac\ne 0 \\
                      |b+[k/h+1/2]a|\leq |a|\leq h \\
                      |d+[k/h+1/2]c|\leq |c|\leq h \\
                      (k/h+b/a)(k/h+d/c)<0 }}
      \sgn\left(\frac{k}{h}+\frac{b}{a}\right)
        (ak+bh)^{\tn}(ck+dh)^{n}.
\end{equation*}
In fact, in the sum $I_{w,n}(h,k)$, each term
\begin{equation*}
      \sgn\left(\frac{k}{h}+\frac{b}{a}\right)
        (ak+bh)^{\tn}(ck+dh)^{n}
\end{equation*}
is equal to zero unless
  $\left(\begin{smallmatrix}a&b\\c&d\end{smallmatrix}\right)\in \Gamma$
satisfies
\begin{equation*}
  |b+[k/h+1/2]a|\leq |a|\leq h,\ \
  |d+[k/h+1/2]c|\leq |c|\leq h.
\end{equation*}
Furthermore we define a function $E_{w,n}:\ZZ^+\times\ZZ\to \CC$
as follows.
\begin{enumerate}
\item
for $n$ odd, $E_{w,n}$ is defined by
  \begin{equation*}
  E_{w,n}(h,k):=I_{w,n}(h,k)
        -\frac{\bar{B}_{n+1}(\frac{k}{h})}{n+1}h^w
        -\frac{\bar{B}_{\tn+1}(\frac{k}{h})}{\tn+1}h^w
        +\frac{w+2}{B_{w+2}}\frac{B_{n+1}}{n+1}\frac{B_{\tn+1}}{\tn+1}h^w;
  \end{equation*}
\item
for $n$ even, $E_{w,n}$ is defined by
  \begin{equation*}
  E_{w,n}(h,k):=I_{w,n}(h,k)
        +\frac{\bar{B}_{n+1}(\frac{k}{h})}{n+1}h^w
        -\frac{\bar{B}_{\tn+1}(\frac{k}{h})}{\tn+1}h^w.
  \end{equation*}
\end{enumerate}
\end{defn}

Since $I_{w,n}(1,0)=0$, we obtain the following directly:
\begin{lem}\label{lem5.1}
If $n$ is odd, we have
\begin{equation*}
  E_{w,n}(1,0)=
    -\frac{B_{n+1}}{n+1}
    -\frac{B_{\tn+1}}{\tn+1}
    +\frac{w+2}{B_{w+2}}\frac{B_{n+1}}{n+1}\frac{B_{\tn+1}}{\tn+1}.
\end{equation*}
\end{lem}

Now we also define a function $S_{w,n}(h,k)$ in $h$ and $k$, which plays a role of
reciprocity function for $E_{w,n}$.
\begin{defn}\label{defn5.2}
Let $n$ be an integer such that $0<n<w$.
We define a polynomial $S_{w,n}$
in $h$ and $k$ as follows.
\begin{enumerate}
\item
for $n$ odd, $S_{w,n}$ is defined by
  \begin{multline*}
    S_{w,n}(h,k):=
        -\frac{B_{n+1}(\frac{k}{h})}{n+1}h^w
        +\frac{B_{n+1}(\frac{h}{k})}{n+1}k^w
        -\frac{B_{\tn+1}(\frac{k}{h})}{\tn+1}h^w
        +\frac{B_{\tn+1}(\frac{h}{k})}{\tn+1}k^w \\
        +\frac{w+2}{B_{w+2}}\frac{B_{n+1}}{n+1}\frac{B_{\tn+1}}{\tn+1}
          (h^w-k^w);
  \end{multline*}
\item
for $n$ even, $E_{w,n}$ is defined by
  \begin{equation*}
    S_{w,n}(h,k):=
        +\frac{B_{n+1}(\frac{k}{h})}{n+1}h^w
        +\frac{B_{n+1}(\frac{h}{k})}{n+1}k^w
        -\frac{B_{\tn+1}(\frac{k}{h})}{\tn+1}h^w
        -\frac{B_{\tn+1}(\frac{h}{k})}{\tn+1}k^w.
  \end{equation*}
\end{enumerate}
\end{defn}
Note that the right hand sides of above equations
are homogeneous polynomials in $h$ and $k$ of degree $w$.
Moreover, we know that the polynomial $S_{w,n}(h,k)$ is even or odd
depending on $n$ is odd or even.
Here are a couple of examples of $S_{w,n}$:
\begin{equation*}
    S_{10,4}(h,k)=-\frac{2h^9k}{35}+\frac{5h^7k^3}{14}-\frac{3h^5k^5}{5}
      +\frac{5h^3k^7}{14}-\frac{2hk^9}{35}.
\end{equation*}
and
\begin{equation*}
    S_{10,5}(h,k)=-\frac{6h^{10}}{691}+\frac{h^8k^2}{6}-\frac{h^6k^4}{2}
    +\frac{h^4k^6}{2}-\frac{h^2k^8}{6}+\frac{6k^{10}}{691}
\end{equation*}

The following is a ``weighted version'' of
\cite[Theorems 1.1, 1.2]{F3}. The proof is similar to
that of the ``non-weighted version'', and we omit it.

\begin{thm}\label{thm5.2}
Let $n$ be an integer such that $0<n<w$.  Then the following assertions hold:
\begin{enumerate}
\item
$E_{w,n}$ is an odd $($resp. even$)$ Dedekind symbol of weight $w$
for $n$ even $($resp. odd$)$.
\item
$E_{w,n}$ has the following reciprocity law:
\begin{equation*}
    E_{w,n}(h,k)-E_{w,n}(k,-h)=S_{w,n}(h,k).
\end{equation*}
\item
Let $E^-$ be an odd Dedekind symbol of weight $w$
whose reciprocity function is a polynomial.
Then  $E^-$ is a linear combination of
 $E_{w,n}\ (0<n< w;\ n \text{\ even})$.
\item
Let $E^+$ be an even Dedekind symbol of weight $w$
whose reciprocity function is a polynomial.
Then  $E^+$ is a linear combination of
 $E_{w,n}\ (0<n< w;\ n \text{\ odd})$,
$F_w$ and $G_w$.
\end{enumerate}
\end{thm}

\section{Dedekind symbols associated with cusp forms of $w\leq 24$}
\label{sect6}

In this section we investigate Dedekind symbols associated
with cusp forms $f\in S_{w+2}$ with $w\leq 24$.
In this case, the dimension of $S_{w+2}$ is at most one.
Then, using Theorem \ref{thm5.2}, we can give explicit formula for
$E_f$.

Let $\ell$ be one of the integer in $\{10,14,16,18,20,24\}$,
then $S_{w+2}$ is one-dimensional for each $\ell$ in this set.
We define $f_{\ell+2}$ to be a unique
normalized eigenform in $S_{\ell+2}$
($a_{f_{\ell+2}}(0)=0$, $a_{f_{\ell+2}}(1)=1$).
Then it is well known that $f_{\ell+2}$ are expressed by
discriminant $\Delta$ and Eisenstein series $Q$, $R$.
Here
\begin{equation*}
  \Delta(z)=e^{2\pi iz}\prod_{n=1}^{\infty}(1-e^{2\pi inz})^{24}
           =\sum_{n=1}^{\infty}\tau(n)e^{2\pi inz},
\end{equation*}
\begin{equation*}
  Q(z)=E_4(z)=1+240\sum_{n=1}^{\infty}\sigma_3(n)e^{2\pi inz}
\end{equation*}
and
\begin{equation*}
  R(z)=E_6(z)=1-504\sum_{n=1}^{\infty}\sigma_5(n)e^{2\pi inz}.
\end{equation*}

Indeed for $\ell=10,14,16,18,20$ or $24$,
$f_{\ell+2}$ are given, respectively, by
$\Delta$, $Q\Delta$, $R\Delta$, $Q^2\Delta$,
$QR\Delta$, $Q^2R\Delta$ (refer to \cite{SD1}).
We use the following notation for the Fourier coefficient of $f_{\ell+2}$:
\begin{equation*}
  f_{\ell+2}(z)=\sum_{n=1}^{\infty}
    \tau_{\ell+2}(n)e^{2\pi inz}.
\end{equation*}

Note that $\tau_{12}$ is nothing but Ramanujan's tau function,
namely $\tau_{12}(n)=\tau(n)$.

We express even Dedekind symbol associated with $f_{\ell+2}$
in terms of $E_{\ell,n}$ which was explicitly given in
Definition \ref{defn5.1}.
\begin{lem}\label{lem6.1}
Let $\ell=10,14,16,18,20$ or $24$.
Then
  \begin{equation*}
  E_{f_{\ell+2}}^+=cE_{\ell,2[(\ell+2)/4]-1}
  \end{equation*}
  where $c$ is a constant.
\end{lem}
\begin{proof}
Let $n_0=2[(\ell+2)/4]-1$, and let $\tn_0=\ell-n_0$.
We know that $0<n_0<\ell$ and
$n_0$ is odd. Then there is $f\in S_{\ell+2}$ such that
$S_f^+(h,k)=S_{\ell,n_0}(h,k)$ (\cite[Theorem 1']{KZ1}).
By Theorem \ref{thm5.2} (2),
the reciprocity polynomial of $E_{\ell,n_0}$ is $S_{\ell,n_0}$
while the reciprocity polynomial of $E_f^+$ is $S_f^+$
by \eqref{eqn2.12}.
Then $S_f^+(h,k)=S_{\ell,n_0}(h,k)$
implies that
  \begin{equation}\label{eqn6.1}
  E_f^+=E_{\ell,n_0}+cG_\ell
  \end{equation}
for some constant $c$ by Lemma \ref{lem2.5}.

Now we will show $c=0$.
We have
\begin{align*}
  E_f^+(1,0)&=\int_{0}^{i\infty}f(z)z^{w}dz
          =S_f^+(1,0)
          =S_{\ell,n_0}(1,0) \\
          &=-\frac{B_{n_0+1}}{n_0+1}
           -\frac{B_{\tn_0+1}}{\tn_0+1}
           +\frac{\ell+2}{B_{\ell+2}}\frac{B_{n_0+1}}{n_0+1}
           \frac{B_{\tn_0+1}}{\tn_0+1} \\
           &\ \ \ \ \text{(by Definition \ref{defn5.2} (1)}).
\end{align*}
On the other hand, by Lemma \ref{lem5.1}, we have
\begin{equation*}
  E_{\ell,n_0}(1,0)
          =-\frac{B_{n_0+1}}{n_0+1}
           -\frac{B_{\tn_0+1}}{\tn_0+1}
           +\frac{\ell+2}{B_{\ell+2}}\frac{B_{n_0+1}}{n_0+1}
           \frac{B_{\tn_0+1}}{\tn_0+1}.
\end{equation*}
These imply $E_f^+(1,0)=E_{\ell,n_0}(1,0)$.
Hence we know that $c=0$ in the equation \eqref{eqn6.1},
and then
  \begin{equation*}
  E_f^+=E_{\ell,n_0}.
  \end{equation*}

Finally, since $\dim S_{\ell+2}=1$, there is a constant
$c'$ such that $f_{\ell+2}=c'f$. This implies
  \begin{equation*}
  E_{f_{\ell+2}}^+=c'E_f^+=c'E_{\ell,n_0}.
  \end{equation*}
This completes the proof.
\end{proof}

\section{Formulae for generalized Ramanujan's tau functions}
\label{sect7}

In this section we calculate $T_mE_{f_{\ell+2}}(h,k)$
and obtain explicit formula for $\tau_{\ell+2}(m)$.

We start with the following lemma.
\begin{lem}\label{lem7.1}
Let $m$ be a positive prime integer,
and let $n$ be odd.
Then $T_mE_{w,n}(h,k)$ is calculated as follows:
  \begin{multline*}
  T_mE_{w,n}(h,k)=
    I_{w,n}(h,mk)+\sum_{b=0}^{m-1}I_{w,n}(mh,k+bh) \\
        \ccsp\ \
          -\frac{1}{n+1}
          \left\{
          \bar{B}_{n+1}(\frac{mk}{h})+m^{\tn}\bar{B}_{n+1}(\frac{k}{h})
          \right\}h^w  \\
        \cccsp\ \ \ \ -\frac{1}{\tn+1}
          \left\{
          \bar{B}_{\tn+1}(\frac{mk}{h})+m^{n}\bar{B}_{\tn+1}(\frac{k}{h})
          \right\}h^w \\
        \ccccsp\ \ +(1+m^{w+1})
        \frac{w+2}{B_{w+2}}\frac{B_{n+1}}{n+1}\frac{B_{\tn+1}}{\tn+1}h^w.
  \end{multline*}
In particular we have
  \begin{multline*}
  T_mE_{w,n}(1,0)=
    \sum_{b=0}^{m-1}I_{w,n}(m,b)
      -\frac{B_{n+1}}{n+1}(1+m^{\tn})
        -\frac{B_{\tn+1}}{\tn+1}(1+m^{n}) \\
      +(1+m^{w+1})
        \frac{w+2}{B_{w+2}}\frac{B_{n+1}}{n+1}\frac{B_{\tn+1}}{\tn+1}.
  \end{multline*}
\end{lem}
\begin{proof}
Applying the formula for
Bernoulli function

\begin{equation*}
  \sum_{b=0}^{c-1}\bar{B}_{n+1}(x+\frac{b}{c})
    =c^{-n}\bar{B}_{n+1}(cx)\ \ \ (c\in\ZZ^+),
\end{equation*}
we have
  \begin{align*}
  T_m&E_{w,n}(h,k) \\
  =&E_{w,n}(h,mk)
    +\sum_{b=0}^{m-1}E_{w,n}(mh,k+bh) \\
  =&
    I_{w,n}(h,mk)
        -\frac{\bar{B}_{n+1}(\frac{mk}{h})}{n+1}h^w
        -\frac{\bar{B}_{\tn+1}(\frac{mk}{h})}{\tn+1}h^w
        +\frac{w+2}{B_{w+2}}\frac{B_{n+1}}{n+1}\frac{B_{\tn+1}}{\tn+1}h^w \\
    &\ +
    \sum_{b=0}^{m-1}
      \Biggl\{
        I_{w,n}(mh,k+bh)
        -\frac{\bar{B}_{n+1}(\frac{k+bh}{mh})}{n+1}(mh)^w
        -\frac{\bar{B}_{\tn+1}(\frac{k+bh}{mh})}{\tn+1}(mh)^w
      \\
      &\ccccsp\csp
      +\frac{w+2}{B_{w+2}}\frac{B_{n+1}}{n+1}\frac{B_{\tn+1}}{\tn+1}(mh)^w
      \Biggr\} \\
  =&
    I_{w,n}(h,mk)+\sum_{b=0}^{m-1}I_{w,n}(mh,k+bh) \\
        &\ccsp\ \
          -\frac{1}{n+1}
          \left\{
          \bar{B}_{n+1}(\frac{mk}{h})+m^{\tn}\bar{B}_{n+1}(\frac{k}{h})
          \right\}h^w  \\
        &\cccsp\ \ \ \ -\frac{1}{\tn+1}
          \left\{
          \bar{B}_{\tn+1}(\frac{mk}{h})+m^{n}\bar{B}_{\tn+1}(\frac{k}{h})
          \right\}h^w \\
        &\ccccsp\ \ +(1+m^{w+1})
        \frac{w+2}{B_{w+2}}\frac{B_{n+1}}{n+1}\frac{B_{\tn+1}}{\tn+1}h^w. \\
  \end{align*}
This completes the proof.
\end{proof}

Next we calculate $T_mE_{w,n}(1,0)/E_{w,n}(1,0)$.

\begin{lem}\label{lem7.2}
Let $m$ be a positive prime integer,
and let $n$ be odd.
Suppose that $E_{w,n}(1,0)\ne 0$.
Then we have
  \begin{multline*}
  \frac{T_mE_{w,n}(1,0)}{E_{w,n}(1,0)}
  =1+m^{w+1}
    +\frac{1}{E_{w,n}(1,0)} \\
    \times\left\{
      \frac{B_{n+1}}{n+1}(m^{\tn}-m^{w+1})
        +\frac{B_{\tn+1}}{\tn+1}(m^{n}-m^{w+1})
      +\sum_{b=0}^{m-1}I_{w,n}(m,b)
    \right\}.
  \end{multline*}
\end{lem}

\begin{proof}
Applying Lemma \ref{lem7.1}, we have
  \begin{align*}
  \frac{T_mE_{w,n}(1,0)}{E_{w,n}(1,0)}
  =&\frac{1}{E_{w,n}(1,0)}
    \Biggl\{
    \sum_{b=0}^{m-1}I_{w,n}(m,b)
      -\frac{B_{n+1}}{n+1}(1+m^{\tn})
        -\frac{B_{\tn+1}}{\tn+1}(1+m^{n})
    \\
    &\cccsp\ccsp\ \ \
      +(1+m^{w+1})
        \frac{w+2}{B_{w+2}}\frac{B_{n+1}}{n+1}\frac{B_{\tn+1}}{\tn+1}
    \Biggr\} \\
  =&\frac{1}{E_{w,n}(1,0)}
    \Biggl\{
    \sum_{b=0}^{m-1}I_{w,n}(m,b)
      -\frac{B_{n+1}}{n+1}(-m^{w+1}+m^{\tn})
    \\
    &\ccsp\csp
      -\frac{B_{\tn+1}}{\tn+1}(-m^{w+1}+m^{n})
      +(1+m^{w+1})E_{w,n}(1,0)
    \Biggr\} \\
  =&1+m^{w+1} \\
    &+\frac{1}{E_{w,n}(1,0)}
    \Biggl\{
      \frac{B_{n+1}}{n+1}(m^{\tn}-m^{w+1})
        +\frac{B_{\tn+1}}{\tn+1}(m^{n}-m^{w+1})
    \\
    &\ccccsp\ccsp
    +\sum_{b=0}^{m-1}I_{w,n}(m,b)
    \Biggr\}. \\
  \end{align*}
This completes the proof.
\end{proof}

Finally we arrive at explicit formulae
for generalized Ramanujan's tau functions:

\begin{thm}\label{thm7.3}
Let $m$ be a positive prime integer.
For $\ell=10,14,16,18,20$ and $24$,
$\tau_{\ell+2}(m)$ are expressed as
\begin{align*}
  \tau_{12}(m)&=1+m^{11}-\frac{691}{6}
  \left\{
    -\frac{1}{126}(m^5-m^{11})+\sum_{b=0}^{m-1}I_{10,5}(m,b)
  \right\}, \\
  \tau_{16}(m)&=1+m^{15}+\frac{3617}{30}
  \left\{
    \frac{1}{120}(m^7-m^{15})+\sum_{b=0}^{m-1}I_{14,7}(m,b)
  \right\}, \\
  \tau_{18}(m)&=1+m^{17}-\frac{43867}{150}
  \Biggl\{
    -\frac{1}{132}(m^7-m^{17})
    +\frac{1}{240}(m^9-m^{17})
  \\
  &\ccccsp\ccsp\csp\ \ \
    +\sum_{b=0}^{m-1}I_{16,7}(m,b)
  \Biggr\}, \\
  \tau_{20}(m)&=1+m^{19}-\frac{174611}{2646}
  \left\{
    -\frac{1}{66}(m^9-m^{19})+\sum_{b=0}^{m-1}I_{18,9}(m,b)
  \right\}, \\
  \tau_{22}(m)&=1+m^{21}+\frac{77683}{1050}
  \Biggl\{
    \frac{691}{32760}(m^9-m^{21})
    -\frac{1}{132}(m^{11}-m^{21})
  \\
  &\ccccsp\ccsp\csp\ \
    +\sum_{b=0}^{m-1}I_{20,9}(m,b)
  \Biggr\}, \\
  \tau_{26}(m)&=1+m^{25}-\frac{657931}{40950}
  \Biggl\{
    -\frac{1}{12}(m^{11}-m^{25})
    +\frac{691}{32760}(m^{13}-m^{25})
  \\
  &\ccccsp\ccsp\csp\
    +\sum_{b=0}^{m-1}I_{24,11}(m,b)
  \Biggr\}. \\
\end{align*}
\end{thm}
\begin{proof}
By Theorem \ref{thm4.1}, we know
\begin{equation*}
  \tau_{\ell+2}(m)=\frac{T_mE_{f_{\ell+2}}^+(h,k)}{E_{f_{\ell+2}}^+(h,k)}.
\end{equation*}
By Lemma \ref{lem6.1}, we know
\begin{equation*}
  E_{f_{\ell+2}}^+=cE_{\ell,2[(\ell+2)/4]-1},
\end{equation*}
and thus we have
\begin{equation*}
  \tau_{\ell+2}(m)=\frac{T_mE_{\ell,2[(\ell+2)/4]-1}(h,k)}
                     {E_{\ell,2[(\ell+2)/4]-1}(h,k)}
\end{equation*}
for $(h,k)$ with ${E_{\ell,2[(\ell+2)/4]-1}(h,k)}\ne 0$.

Direct calculations using Lemma \ref{lem5.1} yield
\begin{equation*}
  E_{\ell,2[(\ell+2)/4]-1}(1,0)
    =-\frac{6}{691},\frac{30}{3617},-\frac{150}{43867},
     -\frac{2646}{174611},\frac{1050}{77683},-\frac{40950}{657931}
\end{equation*}
for $\ell=10,14,16,18,20,24$, respectively.
In particular we know ${E_{\ell,2[(\ell+2)/4]-1}(1,0)}\ne 0$.

Hence we have
\begin{equation*}
  \tau_{\ell+2}(m)=\frac{T_mE_{\ell,2[(\ell+2)/4]-1}(1,0)}
                     {E_{\ell,2[(\ell+2)/4]-1}(1,0)}.
\end{equation*}
Then we use Lemmas \ref{lem7.2},
putting $w=\ell$ and $n=2[(\ell+2)/4]-1$,
to obtain the formulae in Theorem \ref{thm7.3}.
\end{proof}

\begin{rem}\label{rem7.1}
There are other formulae for $\tau_{\ell+2}$
different from ours,
by MacDonald \cite{MA1} $($see also Dyson \cite{D1}$)$
and by Manin \cite{M1}.
\end{rem}

From Theorem \ref{thm7.3} and the fact that $I_{w,n}(h,k)$
is an integer, we rediscover the following congruences
which are obtained by Ramanujan \cite{RA1}, Swinnerton-Dyer \cite{SD1}
and Manin \cite{M1}.
{\allowdisplaybreaks
\begin{align*}
  \tau_{12}(m)&\equiv\sigma_{11}(m)\mod\ 691, \\
  \tau_{16}(m)&\equiv\sigma_{15}(m)\mod\ 3617, \\
  \tau_{18}(m)&\equiv\sigma_{17}(m)\mod\ 43867, \\
  \tau_{20}(m)&\equiv\sigma_{19}(m)\mod\ 283\cdot617, \\
  \tau_{22}(m)&\equiv\sigma_{21}(m)\mod\ 131\cdot593, \\
  \tau_{26}(m)&\equiv\sigma_{25}(m)\mod\ 657931.
\end{align*}
}
(These congruences can be shown first for $m$ prime, and
then for any positive integer due to multiplicative properties
of $\tau_{\ell+2}$ and $\sigma_{\ell+1}$.)

\section{Knopp-Parson-Rosen identities}
\label{sect8}

In this final section we show how naturally
Knopp-Parson-Rosen identities are derived from
Hecke operators on Dedekind symbols.

Let $f\in M_{w+2}$ be a normalized Hecke eigenform, and let $E_f^-$ be
the odd Dedekind symbol associated with $f$.
Then, in Theorem \ref{thm4.1}, we showed the following
identity:
  \begin{equation}\label{eqn8.1}
    \sum_{\substack{ad=n \\ 0<d}}
      \sum_{b(\md d)}E_f^-(dh,ak+bh)
      =a_f(n)E_f^-(h,k)
  \end{equation}
for any $(h,k)\in \ZZ^+\times\ZZ$.

Taking $f$ to be Eisenstein series
$G_{w+2}$ of weight $w+2$, we can see the identity
\eqref{eqn8.1} is nothing
but Knopp-Parson-Rosen identity \cite{K1,P1,PR1}.

Now let us recall generalized Dedekind sums $s_{w+1}(k,h)$ introduced
by Apostol.

\begin{defn}[Apostol \cite{A1,CA1}]\label{defn8.1}
\begin{equation*}
  s_{w+1}(k,h)=\sum_{\mu=0}^{h-1}\frac{\mu}{h}\bar{B}_{w+1}(\frac{\mu k}{h}).
\end{equation*}
\end{defn}

In \cite[Lemma 6.4]{F2}), for Eisenstein series  $G_{w+2}$
defined by
\begin{equation*}
  G_{w+2}(z)=-\frac{B_{w+2}}{2(w+2)}+\sum_{m=1}^{\infty}
      \sum_{n=1}^{\infty}n^{w+1}e^{2\pi imnz},
\end{equation*}
it was shown that
\begin{equation}\label{eqn8.2}
  E_{G_{w+2}}^-(h,k)
    =-\frac{h^{w}}{2(w+1)}s_{w+1}(k,h).
\end{equation}
(Though, in \cite{F2}, \eqref{eqn8.2} was proved for $(p,q)\in V$,
the proof is also valid for any $(h,k)\in\ZZ^+\times\ZZ$.)

Then we have

\begin{thm}[\cite{K1,P1,PR1}]\label{thm8.1}
\begin{equation}\label{eqn8.3}
    \sum_{\substack{ad=n \\ 0<d}}
      \sum_{b(\md d)}d^ws_{w+1}(ak+bh,dh)
    =\sigma_{w+1}(n)s_{w+1}(k,h).
\end{equation}
\end{thm}

\begin{proof}
Since $G_{w+2}$ is a normalized eigenform for $T_n$
with the eigenvalue $\sigma_{w+1}(n)$, we have
  \begin{equation}\label{eqn8.4}
    \sum_{\substack{ad=n \\ 0<d}}
      \sum_{b(\md d)}E_{G_{w+2}}^-(dh,ak+bh)
      =\sigma_{w+1}(n)E_{G_{w+2}}^-(h,k).
  \end{equation}
From \eqref{eqn8.4} and \eqref{eqn8.2}, we obtain \eqref{eqn8.3}.
\end{proof}

\begin{rem}\label{rem8.1}
The reader can easily see that
Theorems
\ref{thm1.1}, \ref{thm1.2}, \ref{thm1.3},
\ref{thm1.4} and \ref{thm1.5}
are special cases of the corresponding Theorems
\ref{thm2.3}, \ref{thm2.4}, \ref{thm3.3},
\ref{thm4.1} and \ref{thm7.3}
for modular forms.
If we restrict these latter theorems to cusp forms,
we rediscover Theorems from \ref{thm1.1} to \ref{thm1.5}.
Hence we omit proofs of those theorems.
\end{rem}

\end{document}